\documentclass[11pt]{article}

\usepackage[english]{babel}
\usepackage{graphicx}
\usepackage{subfig}
\usepackage{amsfonts,amssymb,amsmath,latexsym,amsthm}
\usepackage{multirow}
\usepackage[linesnumbered,ruled,vlined]{algorithm2e}
\usepackage{tikz}
\usetikzlibrary{intersections, calc}

\usepackage{thm-restate}

\SetKwInOut{KwIn}{Input}
\SetKwInOut{KwOut}{Output}

\newtheorem{thm}{Theorem}[section]

\newtheorem{lem}[thm]{Lemma}
\newtheorem{prop}{Proposition}
\newtheorem{obs}[thm]{Observation}

\newtheorem{prob}[thm]{Problem}

\newtheorem{claim}[thm]{Claim}
\newtheorem{fact}{Fact}
\newtheorem{case}{Case}

\theoremstyle{definition}
\newtheorem{defi}{Definition}

\begin{document}

\title{Monochromatic triangle-tilings in dense
graphs without large
independent sets 
\thanks{This work was supported by the National Key Research and Development Program of
China (2023YFA1010203), the National Natural Science Foundation of China (No.12471336), and
the Innovation Program for Quantum Science and Technology (2021ZD0302902).}}
\author{
Xinmin Hou$^{a,b}$, \quad Xiangyang Wang$^{a}$, \quad Zhi Yin$^{a}$\\[2mm]
\small $^{a}$School of Mathematical Sciences,\\
\small University of Science and Technology of China, Hefei 230026, Anhui, China\\
\small $^{b}$ Hefei National Laboratory,\\
\small University of Science and Technology of China, Hefei 230088, Anhui, China\\
\small Email: \texttt{xmhou@ustc.edu.cn (X. Hou)},  
\small \texttt{wangxiangyang@mail.ustc.edu.cn (X. Wang)},\\
\small \texttt{yinzhi@mail.ustc.edu.cn (Z. Yin)}
}

\date{}
\maketitle
\begin{abstract}
Given two graphs $H$ and $G$, an  $H$-tiling is a family of vertex-disjoint copies of $H$ in $G$.
A perfect $H$-tiling covers all vertices of $G$.
The Corrádi–Hajnal theorem (1963) states that an $n$-vertex graph $G$  with minimum degree $\delta(G)\ge 2n/3$ contains a perfect triangle-tiling.
For an $n$-vertex graph $G$ with independence number $\alpha(G)=o(n)$, Balogh, Molla and Sharifzadeh (Random Structures \& Algorithms, 2016) showed that a minimum degree of $(\frac12+o(1))n$ forces a perfect triangle-tiling.
In a 2-edge-colored graph,
Balogh, Freschi, Treglown (European J. Combin. 2026) determined the (asymptotic) minimum degree threshold
for forcing a strong or weak monochromatic triangle-tiling covering a prescribed proportion of the 
vertices: a strong tiling requires all triangles to be in the same color class, while a weak tiling only requires each triangle to be monochromatic. 
In this paper, we combine the conditions from these two lines of work and prove that every $2$-edge-colored $n$-vertex graph $G$ with $\alpha(G)=o(n)$ contains a weak monochromatic triangle-tiling $\Gamma$ of size
\[
|\Gamma|\ge
\begin{cases}
2\delta(G)-n-o(n), & \text{if }\frac12 n\le \delta(G)\le \frac35 n,\\[2mm]
\delta(G)/3-o(n), & \text{if }\delta(G)>\frac35 n.
\end{cases}
\]
Both bounds are asymptotically optimal. We use the degree form regularity lemma in our proof.
\end{abstract}

\section{Introduction}
Dirac-type Problem--which ask which global structures are forced by a given minimum-degree condition--are a central theme in extremal graph theory. Given two graphs $H$ and $G$, an {\it $H$-tiling} (or {\it $H$-packing}) is a family of vertex-disjoint copies of $H$ in $G$. A {\it perfect $H$-tiling} (or {\it $H$-factor}) is an $H$-tiling that covers all vertices of $G$. For triangle packings, the Corrádi–Hajnal theorem\cite{CorradiHajnal1963} states that any $n$-vertex graph $G$ (that is, a graph with $n$ vertices) with minimum degree $\delta(G)\ge 2n/3$ contains a triangle-factor.

In graphs with sublinear independence number $\alpha(G)=o(n)$, the absence of large independent sets forces strong neighbourhood overlap and yields robust local density. A key result in this area is by Balogh, Molla and Sharifzadeh~\cite{BMS16}, who obtained the optimal  threshold for the existence of a triangle-factor under this condition. 

\begin{thm}[Balogh, Molla, Sharifzadeh~\cite{BMS16}]\label{thm:BMS}
For every $\varepsilon>0$, there exist $\alpha>0$ and $n_0$ such that the following holds.
If $n>n_0$ and $G$ is an $n$-vertex graph with $\alpha(G)\le \alpha n$ and
\[
\delta(G)\ge \Big(\tfrac12+\varepsilon\Big)n,
\]
then $G$ contains a triangle-factor.
\end{thm}

Several extensions of Theorem~\ref{thm:BMS} have been established in literature. We list some of them.
\begin{itemize}
\item (Knierim, Su~\cite{KrFactorsLowInd})  Given $\varepsilon>0$ and an integer $r\ge 4$, there exists $\alpha>0$ such that the following holds for all sufficiently large $n$ with $r| n$.
If  $G$ is an $n$-vertex graph satisfying
\[
\delta(G)\ge \left(1-\tfrac{2}{r}+\varepsilon\right)n
\quad\text{and}\quad
\alpha(G)\le \alpha n,
\]
then $G$ contains a $K_r$-factor.

\item (Nenadov, Pehova~\cite{NenadovStructure})
    Let $r > \ell \ge 2$ be integers. For any $\varepsilon > 0$ there exists
$\alpha > 0$ such that if $G$ is a graph with $n$ vertices where $r$ divides
$n$, 
\[
\delta(G) \ge \left(\tfrac{r-\ell}{r-\ell+1} + \varepsilon\right) n
\quad \text{and} \quad
\alpha_\ell(G) \le \alpha n,
\]
then $G$ contains a $K_r$-factor, where, $\alpha_\ell(G)$, the \emph{$\ell$-independence number} of graph $G$ is define as the maximum size of an induced $K_\ell$-free subgraph of $G$.

\item (Chen, Han, Yang~\cite{ChenHanYangCliqueFactors}) 
Let $r\ge5$ and $\varepsilon>0$. For the case of $\ell=r-2$, there exists $\alpha>0$ such that for all
sufficiently large $n$ with $r|n$, every $n$-vertex graph $G$ satisfying
\[
\delta(G)\ge\left(\tfrac12+\varepsilon\right)n
\quad\text{and}\quad
\alpha_{r-2}(G)\le \alpha n
\]
contains a $K_r$-factor.

\item (Chang, Han, Kim, Wang, Yang~\cite{EmbeddingCliqueFactors})     Let integers $r,\ell$ satisfying $r>\ell\ge \tfrac34 r$.
For every $\varepsilon>0$, there exists $\alpha>0$ such that for all sufficiently large
$n$ with $r|n$, every $n$-vertex graph $G$ with
\[
\delta(G)\ge \bigl(\tfrac1{2-\rho_\ell(r-1)}+\varepsilon\bigr)n
\quad\text{and}\quad
\alpha_\ell(G)\le \alpha n
\]
contains a $K_r$-factor. Moreover, this minimum degree condition is asymptotically
best possible. Here,  $\rho_\ell(r)$ is the \emph{Ramsey-Tur\'an density} of $K_r$ under the
$\ell$-independence condition, given by 
\[
\rho_\ell(r)
=
\lim_{\varepsilon\to0}\lim_{n\to\infty}
\frac{\mathrm{RT}_\ell(n,K_r,\varepsilon n)}{\binom{n}{2}},
\]
where $\mathrm{RT}_\ell(n,K_r,\varepsilon n)$ is the maximum number of edges in an
$n$-vertex $K_r$-free graph $G$ with $\alpha_\ell(G)\le \varepsilon n$.


\end{itemize}

Another research direction concerns tilings in edge-colored graphs. In this context,  
it is useful to distinguish between \emph{strong} and \emph{weak} monochromatic tilings.
A {\it strong monochromatic $H$-tiling} requires all copies of $H$ to lie in the same color class, whereas a {\it weak monochromatic $H$-tiling} only insists that each individual copy of $H$ be monochromatic.
Since Burr, Erd\H{o}s and Spencer~\cite{BurrErdosSpencer1975} established that the Ramsey number for $m$ disjoint triangles is 
\[
R(mK_3) = 5m,
\] it follows that any 2-edge-coloring of the complete graph $K_{5m}$ contains a strong monochromatic $K_3$-titling of size $m$.
A natural extension of the result of Burr, Erd\H{o}s and Spencer was given by Moon~\cite{Moon1966DisjointTriangles}, who showed the following: For every integer $m \ge 2$, every $2$-edge-colouring of $K_{3m+2}$ yields a $K_3$-tiling consisting of $m$ monochromatic copies of $K_3$.


In this article, we consider the strong (or weak) monochromatic $K_3$-tilings in dense graphs.
Given graphs $G$ and $H$, define the {\em strong (resp. weak) $H$-tiling number} of graph $G$ as 
$$\nu_s(G, H)=\max_c|\mathcal H|, \, (\text{resp. } \nu_w(G, H)=\max_c|\mathcal H|)$$
where $c$ ranges over all 2-edge-coloring of $G$, and $\mathcal H$  is a strong (resp. weak) monochromatic $H$-tiling in this colored graph.  Balogh, Freschi and Treglown~\cite{RTPDense} established the optimal minimum-degree guarantees for both types of tillings.





\begin{thm}[Balogh, Freschi, Treglown~\cite{RTPDense}]\label{thm:BFT-strong}
There exists $n_0$ such that the following holds for all $n\ge n_0$.
Let $G$ be a $2$-edge-coloured $n$-vertex graph with minimum degree $\delta(G)=\delta$.
Then\\ 
(1) \[\nu_s(G, K_3)\ge
\begin{cases}
\Big\lfloor \dfrac{\delta+1}{5} \Big\rfloor, & \text{if }\delta\ge \tfrac{65}{66}n,\\[2mm]
\Big\lceil \dfrac{5\delta-4n}{2} \Big\rceil, & \text{if }\tfrac45 n \le \delta\le \tfrac56 n.
\end{cases}
\]
Both bounds are best possible.\\
(2) \[\nu_w(G, K_3)\ge 
\begin{cases}
5\delta-4n, & \text{if } \frac{4n}{5}\le \delta \le \frac{5n}{6},\\[2mm]
\Big\lfloor \dfrac{4\delta-3n}{2} \Big\rfloor - o(n), & \text{if } \frac{5n}{6}\le \delta \le \frac{7n}{8},\\[3mm]
\Big\lfloor \dfrac{2\delta-n}{3} \Big\rfloor, & \text{if } \delta \ge \frac{7n}{8}.
\end{cases}
\]
Moreover, these bounds are asymptotically best possible in the first and third ranges, and best possible up to the $o(n)$ term in the middle range.
\end{thm}


This paper focuses on the size of weak monochromatic $K_3$-tilings in dense $2$-edge-colored graphs with sublinear independence number; consequently,
 all subsequent references to ``monochromatic'' tilings should be understood in this weak sense unless stated otherwise.
We determine the optimal minimum-degree threshold for forcing large monochromatic triangle-tilings in $2$-edge-colored graphs under the independence-number condition $\alpha(G)=o(n)$.
The following is our main results. 

\begin{restatable}{thm}{Main}
\label{thm:main}
Let $1\gg\gamma\gg\alpha\gg 1/n$ and let $G$ be an $n$-vertex graph with $\alpha(G)\le \alpha n$.
Then 
\[
\nu_w(G, K_3)\ge
\begin{cases}
2\delta(G)-n-\gamma n, & \text{if }\frac12n \le \delta(G)\le \frac35 n,\\[2mm]
\delta(G)/3-\gamma n, & \text{if }\delta(G) > \frac35 n.
\end{cases}
\]
Moreover, these bounds are asymptotically optimal.
\end{restatable}

\noindent{\bf Remark A:} Let $1\gg \alpha \gg 1/n$. To verify the tightness of bouns in Theorem~\ref{thm:main}, we  construct a $2$-edge-colored graph $(G,\phi)$ (where $\phi$ is a 2-edge coloring of $G$) on $n$ vertices with $\alpha(G)\le \alpha n$ such that every monochromatic $K_3$-tiling has size at most
\[
\begin{cases}
2\delta(G)-n, & \text{if }\frac12n < \delta(G)\le \frac35 n,\\[2mm]
\delta(G)/3, & \text{if }\delta(G) > \frac35 n.
\end{cases}
\]
\noindent{\bf Construction of $(G,\phi)$:}
    Let $\ell=\lceil\frac{n}{n-\delta(G)}\rceil$ and
    \begin{align*}
        V_1\cup V_2\cup \dots\cup V_{\ell-1}\cup V_\ell=V(G)
    \end{align*}
    be a partition such that $|V_1|=|V_2|=\dots=|V_{\ell-1}|=n-\delta(G)$ and $|V_\ell|=n-(\ell-1)(n-\delta(G))$. 
    The edge set of $G$ satisfies: 
    \begin{itemize}
    \item all edges between distinct parts $V_i$ and $V_j$ are present,
    \item each induced subgraph $G[V_i]$ is triangle-free and has independent number $\alpha(G[V_i])\le \alpha n$.
\end{itemize}
 The feasibility of the second item follows from the lower bound for the off-diagonal Ramsey number $R(3,\alpha n)$ given by  Campos, Jenssen, Michelen, Sahasrabudhe~\cite{CamposJenssenMichelenSahasrabudhe2025}:
\[
R(3,k) \ge \left(\frac{1}{3} + o(1)\right)\frac{k^2}{\log k}.
\]
From this, we obtain that
\[
    R(3,\alpha n)=\Omega\!\left(\frac{(\alpha n)^2}{\log(\alpha n)}\right)\gg n,
\]
allowing us to choose $G[V_i]$ accordingly.

Let $\phi$ be a 2-edge-coloring of  $G$ by coloring all edges between $V_1$ and $\overline{V_1}$ red, and all remaining edges blue.  
Then no vertex in $V_1$ can lie in any monochromatic triangle. Hence every monochromatic $K_3$-tiling leaves $V_1$ uncovered, and thus has size at most $(n-|V_1|)/3=\frac{\delta(G)}{3}$.

 Furthermore, we can obtain a more refined bound when $\frac12 n < \delta(G) \le \frac35 n$. In this case, we have $\ell=3$ with part sizes $|V_1|=|V_2|=n-\delta(G)$ and $|V_3|=2\delta(G)-n$. Since $\delta(G) \le \frac35n$. we have $2|V_3|=4\delta(G)-2n \le n-\delta(G)=|V_2|$. 
    As before, no vertex in $V_1$ can be covered by any monochromatic triangle. Consequently, every triangle in a maximum monochromatic $K_3$-tiling must lie entirely within $G[V_2\cup V_3]$.
    Since $G[V_2]$ is triangle-free, every monochromatic $K_3$ uses at least one vertex from $V_3$. Hence the size of any monochromatic $K_3$-tiling is at most $|V_3|=2\delta(G)-n$.

\begin{figure}[htbp]
    \centering
    \qquad\subfloat[$\delta(G) \ge \frac{3}{5}n$.]{
    \begin{tikzpicture}[edge/.style = {draw, fill = #1!50},]
    \draw[edge = red] (0,0) ellipse (0.8 and 1.6) node[above = 1.6cm,] {$V_1$};
    \draw[edge = red] (3,0) ellipse (1 and 2);
    \path (0.8,0) arc (0:45:0.8 and 1.6) coordinate (a){}
        (0.8,0) arc (0:-45:0.8 and 1.6) coordinate (b){}
        (2,0) arc (180:135:1 and 2) coordinate (c) {}
        (2,0) arc (180:225:1 and 2) coordinate (d) {};
        \fill[edge = blue] (a) -- (c) arc(135:225:1 and 2) (d) -- (b) arc(-45:45:0.8 and 1.6) (a)
        (0,-2.4);

    \path (4,0) arc (0:30:1 and 2) coordinate (e1){}
    (2,0) arc (180:150:1 and 2) coordinate (e2){};
    \path (4,0) arc (0:-30:1 and 2) coordinate (e3){}
    (2,0) arc (180:210:1 and 2) coordinate (e4){};
    \draw (e1) -- (e2);
    \draw (4,0) -- (2,0);
    \draw (e3) -- (e4);
    \node[node font = \footnotesize] at (3,1.4) {$V_2$};
    \node[node font = \footnotesize] at (3,0.5) {$V_3$};
    \node[node font = \footnotesize] at (3,-0.5) {$\dotsb$};
    \node[node font = \footnotesize] at (3,-1.4) {$V_{l}$};
    \node[node font = \footnotesize] at (0,0) {$n - \delta(G)$};
\end{tikzpicture}
}\hfill
    \subfloat[$\frac{n}{2}\le \delta(G)\le \frac{3}{5}n$.]{
    \begin{tikzpicture}[edge/.style = {draw, fill = #1!50}]
    \draw[edge = red] (2,2.2) circle(0.8) node[above = 0.8cm,] {$V_3$};
    \draw[edge = red] (.4,0) circle(1) node[below = 1cm,] {$V_1$};
    \draw[edge = red] (3.6,0) circle(1) node[below = 1cm,] {$V_2$};
    \fill[edge = blue] (2,2.2) ++(200:0.8) -- ($(.4,0) + (80:1)$) arc (80:40:1) -- ($(2,2.2) + (250:.8)$) arc (250:200:.8);
    \fill[edge = red] (2,2.2) ++(-10:0.8) -- ($(3.6,0) + (100:1)$) arc (100:150:1) -- ($(2,2.2) + (-80:.8)$) arc (-80:-10:.8);
    \fill[edge = blue] (.4,0) ++(20:1) -- ($(3.6,0) + (160:1)$) arc (160:200:1) -- ($(.4,0) + (-20:1)$) arc (-20:20:1);
    \node[node font = \footnotesize] at (.4,0) {$n - \delta(G)$};
    \node[node font = \footnotesize] at (3.6,0) {$n - \delta(G)$};
    \node[node font = \footnotesize] at (2,2.2) {$2\delta(G) - n$};
\end{tikzpicture}
}\qquad
    \caption{Construction of $(G,\phi)$.}
\end{figure}
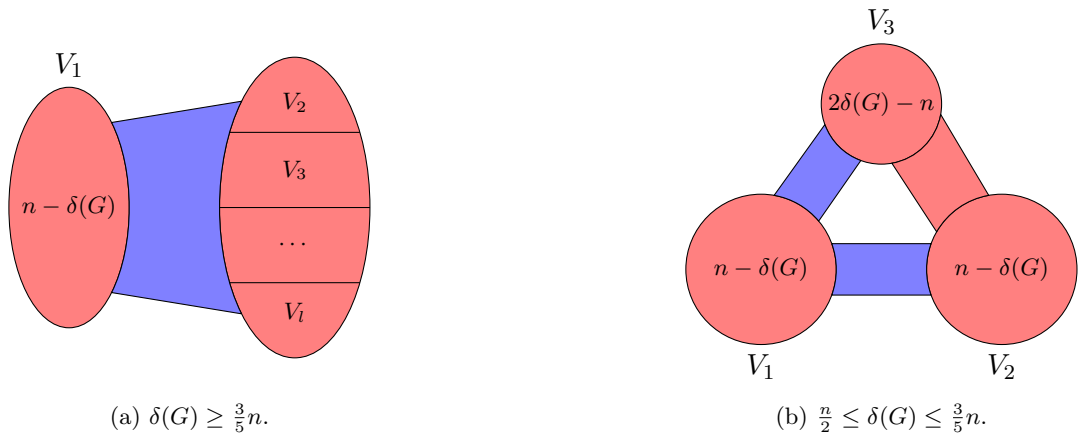

The rest of this article is arranged as follows. Section~2.1 introduces notation and preliminary tools,
and Section~2.2 presents some useful lemmas, including the key Lemma~\ref{lem:k3-tiling-in-F}.
The proof of our main theorem, Theorem~\ref{thm:main}, is given in Section~3, while the proof of  Lemma~\ref{lem:k3-tiling-in-F} is provided in  Section~4.

\section{Preliminary}
\subsection{Notation and definitions}

We now introduce some notation.  
For a graph $G$ and a vertex set $W\subseteq V(G)$, let  
\[
E_G(W)=\{uv\in E(G): u,v\in W\}
\]
denote the set of edges induced by $W$.
Write $(G,\phi)$ for a $2$-edge-colored graph $G$ with 2-edge coloring $\phi$. Define two spanning subgraphs $R$ and $B$ by
\[
E(R)=\{e\in E(G):\phi(e)=1\},
\text{ and }  
E(B)=\{e\in E(G):\phi(e)=2\}, \text{ respectively}.
\]
For each $u\in V(G)$, write
\[
N_R(u)=\{v: uv\in E(R)\},\text{ and } 
N_B(u)=\{v: uv\in E(B)\}.
\]
Clearly, $E(R)\cap E(B)=\emptyset$ and $E(R)\cup E(B)=E(G)$.

Throughout the paper we use the following notation.
\begin{itemize}
    \item $F_2$ denotes the graph consisting of two triangles intersecting in exactly one vertex.
 
    \item $\delta(H,n)$ denotes the smallest integer $k$ such that every $n$-vertex graph $G$ 
    with $\delta(G)\ge k$ and $|H|\mid n$ contains a perfect $H$-tiling.
\end{itemize}

\medskip

Following Kühn and Osthus~\cite{KO}, we recall the parameters that determine 
the minimum-degree threshold for perfect $H$-tilings.

\begin{defi}[Critical chromatic number]
Let $H$ be a graph, and let $\sigma(H)$ denote the minimum size of the smallest color class
over all proper $\chi(H)$-colorings of $H$.  
The \emph{critical chromatic number} of $H$ is
\[
\chi_{cr}(H)
    := \frac{(\chi(H)-1)\,|H|}{\,|H|-\sigma(H)\,}.
\]
\end{defi}

To deal with arithmetic obstructions we use the “highest common factor’’ parameters
from~\cite{KO}.  
Let $\chi(H)=r$. 
For a proper $r$-coloring $c$ with color classes $X_1,\dots,X_r$, define
\[
D(c):=\left\{|\,|X_i|-|X_j|\,|:1\le i,j\le r\right\},
\text{ and } 
D_\chi(H):=\bigcup_{c} D(c).
\]
If $D_\chi(H)\setminus\{0\}\ne\emptyset$, set
\[
\mathrm{hcf}_\chi(H):=\gcd\bigl(D_\chi(H)\setminus\{0\}\bigr),
\]
otherwise, set $\mathrm{hcf}_\chi(H):=\infty$.  
Similarly, let $C_1,\dots,C_s$ be the vertex sets of the connected components of $H$, and define  
\[
D_c(H):=\bigl\{|\,|C_i|-|C_j|\,|:1\le i,j\le s\bigr\}.
\]
We then define $\mathrm{hcf}_c(H)$ as the greatest common divisor  of $D_c(H)\setminus\{0\}$ if this set is nonempty, 
and as $\infty$ otherwise.  Finally, define
\[
\mathrm{hcf}(H):=\gcd\bigl(\mathrm{hcf}_\chi(H),\mathrm{hcf}_c(H)\bigr),
\]
with the conventions $\gcd(\infty,t)=t$ and $\gcd(\infty,\infty)=\infty$.

\begin{defi}[$\chi^{*}(H)$]
For any graph $H$, define
\[
\chi^{*}(H)=
    \begin{cases}
        \chi_{cr}(H), & \text{if }\mathrm{hcf}(H)=1,\\[4pt]
        \chi(H),      & \text{otherwise}.
    \end{cases}
\]
\end{defi}

\begin{thm}[Kühn--Osthus~\cite{KO}]\label{thm:KO}
For every graph $H$ there exists a constant $C_H$ such that
\[
\Bigl(1-\frac1{\chi^{*}(H)}\Bigr)n - 1
\;\le\;
\delta(H,n)
\;\le\;
\Bigl(1-\frac1{\chi^{*}(H)}\Bigr)n + C_H.
\]
\end{thm}

Recall that $F_2$ consists of two triangles intersecting in one vertex.
Then  every proper coloring of $F_2$ has color-class sizes $(1,2,2)$.  Hence $\mathrm{hcf}_\chi(F_2)=1$ and $\chi_{cr}(F_2)=\frac{5}{2}$. Since $F_2$ is connected, 
we have $\mathrm{hcf}_c(F_2)=\infty$. Consequently, $\mathrm{hcf}(F_2)=1$. Therefore $\chi^{*}(F_2)=\chi_{cr}(F_2)=\frac{5}{2}$. 
As a corollary of Theorem~\ref{thm:KO}, we have the following observation.

\begin{obs}\label{obs:F-parameters}
$\delta(F_2,n) \le \frac35n+C_{F_2}$.
\end{obs}

\medskip

\begin{defi}[$\varepsilon$-regularity]
Let $A,B\subseteq V(G)$ be disjoint.
The pair $(A,B)$ is \emph{$\varepsilon$-regular} if for all $X\subseteq A$ and $Y\subseteq B$ with
$|X|\ge \varepsilon |A|$ and $|Y|\ge \varepsilon |B|$,
\[
\bigl|d(X,Y)-d(A,B)\bigr|\le \varepsilon.
\]
where $d(X,Y)=\frac{e(X,Y)}{|X||Y|}$.
\end{defi}
The next two properties were given in \cite{KomlosSimonovits1996Regularity} by Komlós and Simonovits, where the first one is the most important property of regular pairs and the second one says that subgraphs of a regular pair are regular.

\begin{fact}[Fact 1.3 of \cite{KomlosSimonovits1996Regularity}]\label{deg:fact} 
Let $(A,B)$ be an $\varepsilon$-regular pair with density $d$.
Then for any $Y \subseteq B$ with $|Y| > \varepsilon |B|$ we have
\[
\lvert\{x \in A : \lvert N(x)\cap Y\rvert \le (d-\varepsilon)|Y|\}\rvert \le \varepsilon |A|.
\]
More generally, if we fix a set $Y \subseteq B$, and $\ell$ vertices
$x_i \in A$, then ``typically'' they have at least the expected
$d^{\ell}|Y|$ neighbours in common.
\end{fact}

\begin{fact}[Slicing Lemma, Fact 1.5 of \cite{KomlosSimonovits1996Regularity}]\label{slicing:fact}
Let $(A,B)$ be an $\varepsilon$-regular pair with density $d$, and, for some
$\alpha > \varepsilon$, let $A' \subseteq A$ with $|A'| \ge \alpha |A|$ and
$B' \subseteq B$ with $|B'| \ge \alpha |B|$. Then $(A',B')$ is an
$\varepsilon'$-regular pair with
\[
\varepsilon' = \max\{\varepsilon/\alpha,\, 2\varepsilon\},
\]
and for its density $d'$ we have
\[
|d' - d| < \varepsilon.
\]
\end{fact}

\begin{lem}[Regularity Lemma, degree form \cite{KomlosSimonovits1996Regularity}]\label{reg:lem} 
For every $\varepsilon>0$ and every $k_0\in\mathbb{N}$ there exists $M=M(\varepsilon,k_0)$ such that the following holds.
For any graph $G$ on $n$ vertices and any $\beta\in[0,1]$,  
there exist a partition
\[
V(G)=V_0\cup V_1\cup\cdots\cup V_k
\]
and a spanning subgraph $G'\subseteq G$ such that:
\begin{itemize}
    \item[\rm(a)] $k_0 \le k \le M$;
    \item[\rm(b)] $|V_0|\le \varepsilon n$;
    \item[\rm(c)] $|V_1|=\cdots=|V_k|=m\le \varepsilon n$;
    \item[\rm(d)] $d_{G'}(v)> d_G(v)-(\beta+\varepsilon)n$ for all $v\in V(G)$;
    \item[\rm(e)] each $V_i$ is independent in $G'$;
    \item[\rm(f)] each pair $(V_i,V_j)$ is $\varepsilon$-regular in $G'$ with density either $0$ or at least $\beta$.
\end{itemize}
\end{lem}

\begin{defi}[Reduced graph]\label{rdg:def}
Given $G$, $\beta$, $\varepsilon$ and the corresponding subgraph $G'$ from Lemma~\ref{reg:lem},
the \emph{reduced graph} $R_{\beta,\varepsilon}$ is the graph with vertex set $\{1,\dots,k\}$ in which
$ij$ is an edge if and only if $d_{G'}(V_i,V_j)\ge \beta$.
\end{defi}

\begin{prop}\label{deg:prop}
If $R$ is the reduced graph of $G$ obtained using Lemma~\ref{reg:lem} with parameters $\beta,\varepsilon$, then
\[
\delta(R)\;\ge\;
\left(\frac{\delta(G)}{n}-(\beta+\varepsilon)\right)\,|V(R)|.
\]
\end{prop}

\begin{proof}
Applying  Lemma~\ref{reg:lem} to $G$ and let $G'$ be the resulting subgraph.
Since $d_{G'}(v) > d_G(v)-(\beta+\varepsilon)n$, we have
\[
|E_{G'}(V_j,\overline{V_j})|
  = \sum_{x\in V_j} d_{G'}(x)
  \ge m\bigl(\delta(G)-(\beta+\varepsilon)n\bigr).
\]
On the other hand, by the definition of $R$,
\[
|E_{G'}(V_j,\overline{V_j})|
   \le m^2\,\delta(R),
\]
where $j$ is a cluster of minimum degree in $R$.
Combining the two bounds yields
\[
m^2\delta(R)
   \ge m\bigl(\delta(G)-(\beta+\varepsilon)n\bigr).
\]
Since $n/m\ge |V(R)|$, the claimed inequality follows.
\end{proof}

\subsection{Lemmas}
The following is a simple observation.
\begin{obs}\label{K3:free}
Let $(G,\phi)$ be a $2$-edge-colored graph.  
If  $E_R(N_R(u))\neq\emptyset$ (or $E_B(N_B(u))\neq\emptyset$) for some $u\in V(G)$,  
then $(G,\phi)$ contains a monochromatic $K_3$ containing $u$.
Moreover, if for some $u,v\in V(G)$,
\[
|N_R(u)\cap N_B(v)|\ge \alpha(G)+1,
\]
then $(G,\phi)$ also contains a monochromatic copy of $K_3$.
\end{obs}

\begin{lem}\label{lem:F-tiling}
Let $C_{F_2}$ be the constant given by Observation~\ref{obs:F-parameters} for $F_2$.
Then there exists another constant $C=C(F_2)$ such that the following holds.
Let $R$ be a graph on $k$ vertices with $\delta(R)\le \frac{3}{5}k$.  
Then $R$ contains an $F_2$-tiling $\mathcal{T}$ with
\begin{align}\label{F-tiling:eq}
    |\mathcal{T}| \;\ge\; 2\delta(R) - k - C.
\end{align}
\end{lem}

\begin{proof}
Let $\delta=\delta(R)$.  
Since the inequality~\eqref{F-tiling:eq} holds trivially when $\delta \le \frac12k$, we assume that $\delta > \frac12 k$. 
Choose a constant $C>0$ so that (1) $C \ge \frac{5}{2}C_{F_2} + 10$; (2) $\frac32 k - \frac52 \delta + C \in \mathbb{N}$;
and   (3) $5$ divides \(\tfrac52 k - \tfrac52 \delta + C\).
Define an auxiliary graph $R'$ by
\[
V(R') = V(R)\cup W,
\text{ where }
|W| = \frac32 k - \frac52 \delta + C,
\]
and
\[
E(R')
 = E(R)\cup\{ab : a\in V(R),\, b\in W\}.
\]
Thus every $w\in W$ has degree $k$, and every $v\in V(R)$ has degree
\[
d_{R'}(v) = d_R(v) + |W|
           \ge \delta + |W|
           = \frac32 k - \frac32\delta + C.
\]
Hence according to $\delta(R) \in (\frac12k,\frac35k)$, we have
\[
\delta(R') \ge \min\left\{k,\frac32 k - \frac32\delta + C\right\} = \frac32 k - \frac32\delta + C.
\]
Moreover,
\[
|V(R')|
 = k + |W|
 = \frac52 k - \frac52\delta + C.
\]
A direct calculation yields
\[
\delta(R')
  \ge \frac35 |V(R')| + \frac25 C.
\]
Since $C\ge \frac52 C_{F_2}$, we have
\[
\delta(R') \;\ge\; \frac35 |V(R')| + C_{F_2}.
\]

By the choice of $C$, the number $|V(R')|$ is divisible by $5$.  
According to Observation~\ref{obs:F-parameters}, the graph $R'$ contains a
perfect $F$-tiling $\mathcal{T}'$.
Since $W$ is an independent set of $R'$ and $\alpha(F_2)=2$, we conclude that $\mathcal{T}'$ consists of three types of copies of $F_2$:
\begin{itemize}
    \item $s$ copies using exactly two vertices of $W$;
    \item $t$ copies using exactly one vertex of $W$;
    \item $\ell$ copies contained entirely in $R$.
\end{itemize}
Since the first type uses exactly $3$ vertices of $R$,
and the second one uses exactly $4$ vertices of $R$, 
we obtain that
$$\begin{cases}\label{stl:eq}
    2s + t = |W| = \frac32 k - \frac52\delta + C, \\ 3s + 4t + 5\ell = |V(R)| = k.
\end{cases}$$
Solving this linear equation system, we obtain that
\[
\ell - s = 2\delta - k - \frac{4}{5}C.
\]
In particular,
\[
\ell \;\ge\; 2\delta - k - C.
\]
Therefore, the $\ell$ copies of $F_2$ that lie entirely inside $R$ form an $F_2$-tiling
$\mathcal{T}$ in $R$ with
\[
|\mathcal{T}| = \ell \;\ge\; 2\delta(R) - k - C.
\]
Since $C$ depends only on $F_2$, the lemma holds.
\end{proof}


\begin{restatable}{lem}{mainlem}\label{lem:k3-tiling-in-F}
Given $1 \gg \beta \gg \varepsilon \gg \frac{m}{n} \gg \alpha \gg \frac{1}{n}$.
Let $V(G) = V_1 \cup V_2 \cup V_3 \cup V_4 \cup V_5$ with $|V_i| = m$ for all $i \in [5]$,  
and assume $\alpha(G) < \alpha n$.  
Suppose all pairs 
\[
(V_1,V_2),\ (V_2,V_3),\ (V_1,V_3),\ (V_1,V_4),\ (V_4,V_5),\text { and } (V_1,V_5)
\]
are $\varepsilon$–regular, and moreover $d(V_i,V_j) > \beta$.
Let $\phi$ be a $2$–edge–coloring of $G$.  
Then there exists a monochromatic $K_3$–tiling $\Gamma$ in $G$ such that
\[
|\Gamma| \ge (1-\sqrt{\varepsilon})m.
\]
\end{restatable}

\section{Proof of Theorem~\ref{thm:main}}
We restate our main theorem here.
\Main*

\begin{proof}
Let $\gamma \gg \beta \gg \varepsilon \gg \frac{1}{k_0} \gg \alpha \gg \frac{1}{n}$.
By Lemma~\ref{reg:lem} and Definition~\ref{rdg:def}, there exists a constant $M=M(\varepsilon,k_0)$ and a spanning subgraph $G' \subseteq G$
and a partition
\[
V(G) =V_0 \cup V_1 \cup \cdots \cup V_k, \quad k_0 \le k \le M
\]
satisfying all six items in Lemma~\ref{reg:lem}. 
Note that we can choose $\alpha$ small enough such that $\alpha\ll \frac{m}{n}\le \frac1{k_0}\ll\varepsilon$. Then  for $i \in [k]$, $
|V_i|:=m\le \frac{n}{k_0}.
$

Let $R = R_{\beta,\varepsilon}$ be the reduced graph.
From Proposition~\ref{deg:prop}, it follows
\[
\delta(R)
\;\ge\;
\left(\frac{\delta(G)}{n}-(\beta+\varepsilon)\right)\,k.
\]

\begin{case}{$\delta(G) \le \frac{3}{5}n$}.
\end{case}
Applying Lemma~\ref{lem:F-tiling} to $R$,  
there exists a constant $C=C(F_2)$ and
the reduced graph $R$ admits an $F_2$–tiling of size at least $2\delta(R) - k - C$. 
Clearly, each $F_2$-copy in $R$ corresponds to a pair of $5$ vertex sets in $G$, denoted as $V_1,V_2,V_3,V_4,V_5$, where $$(V_1,V_2),(V_1,V_3),(V_2,V_3),(V_1,V_4),(V_1,V_5),(V_4,V_5)$$ are $\varepsilon$-regular, and have density no less than $\beta$.
Let $\phi$ be a 2-edge-coloring of $G$.
Applying Lemma~\ref{lem:k3-tiling-in-F} to the graph $H:=G\left[\ \bigcup_{i=1}^{5}V_i\ \right],$
we obtain that $(H,\phi|_{H})$ admits a monochromatic $K_3$-tiling with size at least $(1-\sqrt{\varepsilon})m$.

Recall that $R$ admits an $F_2$-tiling of size at least $2\delta(R) - k - C$. Thus $(G,\phi)$ contains a monochromatic $K_3$-tiling $\Gamma$ with size at least
\begin{align*}
    |\Gamma|&\ge (1-\sqrt{\varepsilon})m \cdot (2\delta(R) - k - C)\\
    &\ge (1-\sqrt{\varepsilon})(n-\varepsilon n)\left(\frac{2\delta(G)-n}{n}-2(\beta+\varepsilon)-\frac{C}{k}\right)\\
    &\ge (1-2\sqrt{\varepsilon})(2\delta(G)-n-3\beta n)\\
    &\ge 2\delta(G)-n-\frac\gamma2 n,
\end{align*}
where the second inequality holds since
$$km = \sum_{1\le i \le k}|V_i| =n -|V_0| \ge n-\varepsilon n,$$
and last two inequalities derived from the choice of parameters: 
$$\gamma \gg \beta \gg \varepsilon \gg\frac{1}{k_0} \ge \frac1k.$$

\begin{case}
$\delta(G) > \frac{3}{5}n$.
\end{case}
We claim that $(G,\phi)$ contains a monochromatic $K_3$-copy  in this case.  Since
$\chi(F_2)=3$ and
$$\delta(R) \ge (\frac35-\beta-\varepsilon)k > (\frac12+o(1)) k= (1-\frac1{\chi(F_2)-1}+o(1))k,$$
$R$ contains a $F_2$-copy by the Erd\H{o}s-Stone Theorem~\cite{ES46}. Then we can find a monochromatic $K_3$-copy in $(G,\phi)$ by the arguments as in the first paragraph of the above case.

We delete monochromatic $K_3$-copies one by one from $G$, until the resulting subgraph $G'$ satisfying 
\begin{align}\label{delete1:eq}
    \frac35|V(G')|-3\le \delta(G') < \frac35|V(G')|.
\end{align}
Hence  $(G[V(G)\setminus V(G')],\phi)$ contains a monochromatic $K_3$-tiling $\Gamma'$ of size $\frac13(n-|V(G')|)$.
According to the definition of $G'$,  we have $$\delta(G') \ge \delta(G)-(n-|V(G')|).$$ 
Combining with \eqref{delete1:eq}   
we obtain
\begin{align}\label{delete2:eq}
        |V(G')|<\frac52\left(n-\delta(G)\right).
    \end{align}


Recall that $\delta(G') <\frac35 |V(G')|$. Applying Case 1 to $(G',\phi)$ yields a monochromatic $K_3$-tiling $\Gamma''$ of size at least 
$$
2\delta(G')-|V(G')|-\frac\gamma2 |V(G')|.
$$

Now we take the union of these two disjoint tilings, i.e. let $\Gamma=\Gamma'\cup \Gamma''$. Its size 
\begin{align*}
|\Gamma| &\ge \frac13(n-|V(G')|)+2\delta(G')-|V(G')|-\frac\gamma2 |V(G')|\\
&\ge \frac13 n+2\delta(G')-\frac43 |V(G')|-\frac\gamma2 n\\
&\ge \frac13 n+\frac65 |V(G')|-6-\frac43 |V(G')|-\frac\gamma2 n\\
&=  \frac13 n-\frac2{15} |V(G')|-6-\frac\gamma2 n\\
&> \frac13 n-\frac13(n-\delta(G))-6-\frac\gamma2 n\\
&\ge \frac13\delta(G)-\gamma n.
\end{align*}

This completes the proof.
\end{proof}

\section{Proof of Lemma~\ref{lem:k3-tiling-in-F}}
\begin{lem}[Dominating Lemma]\label{lem:dominating}
Given $d > 2\varepsilon>0$ and $t=t(d,\varepsilon)$.
Let $G$ be a bipartite graph with bipartite parts $V_1$ and $V_2$.  
Assume that $(V_1, V_2)$ is $\varepsilon$–regular with density at least $d$, and has sufficiently large size.  
Then there exist vertices $u_1, u_2, \ldots, u_t \in V_1$ such that
\[
\left| \bigcup_{i=1}^{t} N(u_i) \right| \ge (1-\varepsilon)\,|V_2|.
\]
\end{lem}

\begin{proof}
Fix $0<\varepsilon<d/2$ and set $\gamma := d-2\varepsilon>0$.  
Let $t=t(d,\varepsilon)$ be the smallest integer such that
\[
(1-\gamma)^t < \varepsilon.
\]

We construct vertices $u_1,u_2,\dots,u_t\in V_1$ inductively.  
Put $V_2^{(0)}:=V_2$.  
At step $0 \le i\le t$, assume that $u_1,\dots,u_i$ have been chosen (pairwise distinct) and define
\[
V_2^{(i)} := V_2 \setminus \bigcup_{j=1}^i N(u_j).
\]

If $|V_2^{(i)}|\le \varepsilon|V_2|$, then
\[
\Bigl|\bigcup_{j=1}^i N(u_j)\Bigr|
 = |V_2|-|V_2^{(i)}|
 \ge (1-\varepsilon)|V_2|.
\]
Since adding arbitrary vertices
$u_{i+1},\dots,u_t\in V_1$ does not decrease the union, and thus
\[
\Bigl|\bigcup_{j=1}^t N(u_j)\Bigr|\ge (1-\varepsilon)|V_2|.
\]
So from now on we may assume $|V_2^{(i)}|>\varepsilon|V_2|$.

To choose $u_{i+1}$, consider the set of ``bad'' vertices
\[
B := \{ v\in V_1 : |N(v)\cap V_2^{(i)}| < \gamma |V_2^{(i)}| \}.
\]
We claim that $|B|<\varepsilon|V_1|$. If not, then since $|V_2^{(i)}|>\varepsilon|V_2|$, applying
$(d,\varepsilon)$--regularity to $B\subseteq V_1$ and $V_2^{(i)}\subseteq V_2$ yields
\[
d(B,V_2^{(i)}) \ge d(V_1,V_2)-\varepsilon \ge d-\varepsilon.
\]
But by the definition of $B$,
\[
d(B,V_2^{(i)})
 = \frac{e(B,V_2^{(i)})}{|B||V_2^{(i)}|}
 < \gamma = d-2\varepsilon,
\]
a contradiction.

Define the ``good'' set
\[
D := V_1\setminus B.
\]
Then
\[
|D|>(1-\varepsilon)|V_1|
\quad\text{and}\quad
|N(v)\cap V_2^{(i)}|\ge \gamma |V_2^{(i)}|
\quad\text{for all } v\in D.
\]
Since only $i\le t-1$ vertices have been chosen so far, and $|V_1|$ is large,
we may select $u_{i+1}\in D\setminus\{u_1,\dots,u_i\}$.
This vertex satisfies
\[
|N(u_{i+1})\cap V_2^{(i)}| \ge \gamma |V_2^{(i)}|.
\]
Define
\[
V_2^{(i+1)} := V_2^{(i)}\setminus N(u_{i+1}).
\]
Then
\[
|V_2^{(i+1)}|
 = |V_2^{(i)}| - |N(u_{i+1})\cap V_2^{(i)}|
 \le (1-\gamma)\,|V_2^{(i)}|.
\]

As long as $|V_2^{(i)}|>\varepsilon|V_2|$ we may continue this construction.
Thus either $|V_2^{(i)}|\le\varepsilon|V_2|$ for some $i\le t$, or we reach $i=t$.
In the first case, we already showed that
\[
\Bigl|\bigcup_{j=1}^t N(u_j)\Bigr|\ge (1-\varepsilon)|V_2|.
\]
In the second case, induction gives
\[
|V_2^{(i)}|\le (1-\gamma)^i |V_2| \qquad\text{for all } i\le t,
\]
and hence
\[
|V_2^{(t)}| \le (1-\gamma)^t |V_2| < \varepsilon|V_2|.
\]
Thus
\[
\left|\bigcup_{j=1}^t N(u_j)\right|
 = |V_2| - |V_2^{(t)}|
 \ge (1-\varepsilon)|V_2|.
\]

This completes the proof.
\end{proof}

\begin{lem}\label{lem:three-part}
Let $\beta \gg \varepsilon \gg \frac{m}{n} \gg \alpha \gg \frac{1}{n}$.
Let $G$ be a graph whose vertex set admits a partition $V(G)=P\cup Q\cup S$ with \[
n \ge |P|,|Q|,|S| \;\ge\; \sqrt\varepsilon m,
\text{ and }
\alpha(G) < \alpha n.
\] Suppose that each of the bipartite subgraphs induced by the pairs $(P,Q)$, $(Q,S)$, $(S,P)$ is $\sqrt\varepsilon$–regular  
with densities at least $\beta - \varepsilon$.
Then, for any $2$–edge–coloring $\varphi$ of $G$,  
the colored graph $(G,\varphi)$ contains a monochromatic triangle $K$ satisfying
\[
|V(K) \cap S| \le 1,
\text{ and } 
|V(K)\cap P|, |V(K)\cap Q| \le 2.
\]
\end{lem}

\begin{proof}
Assume, for contradiction, that no monochromatic triangle $K$ has a vertex distribution of any of the following types:
\begin{enumerate}
    \item[\rm(1)] $|V(K)\cap P| = 1$ and $|V(K)\cap Q| = 2$;
    \item[\rm(2)] $|V(K)\cap P| = 2$ and $|V(K)\cap Q| = 1$;
   \item[\rm(3)] $|V(K)\cap P| = 2$ and $|V(K)\cap S| = 1$;
    \item[\rm(4)] $|V(K)\cap Q| = 2$ and $|V(K)\cap S| = 1$;
    \item[\rm(5)] $|V(K)\cap P| = |V(K)\cap Q| = |V(K)\cap S| = 1$.
\end{enumerate}
\begin{claim}\label{claim:two-part}
Given two distinct sets $X,Y \in \{P,Q,S\}$. 
If there is no monochromatic triangle $K$ with $|V(K)\cap X|=2$ and $|V(K)\cap Y|=1$,
then there exists an integer $h = h(\beta,\varepsilon)$  
and some subsets
\[
X_1^{r},\ldots,X_h^{r},\;
X_1^{b},\ldots,X_h^{b}\subseteq X
\text{ with }
X^{r}:=\bigcup_{i=1}^h X_i^{r} \text{ and }
X^{b}:=\bigcup_{i=1}^h X_i^{b},
\]
that satisfy the following:

\begin{enumerate}
    \item[\rm(1)] 
    $|X_i^{r}\cap X_j^{b}| < \alpha n$ for all $i\neq j$,  
    and
    \(
    |X^{r}\cup X^{b}| > (1-\varepsilon)|X|.
    \)

    \item[\rm(2)] 
    $E(X_i^{r}) \cap E(R) = E(X_i^{b}) \cap E(B) = \emptyset$ for all $1\le i\le h$.

    \item[\rm(3)] 
    \(
    |E(X^{r}, Y) \cap E(R)|,
    |E(X^{b}, Y) \cap E(B)| \le h \alpha n^2.
    \)
\end{enumerate}
 
\end{claim}

\begin{proof}
By Lemma~\ref{lem:dominating}, there exist $h=h(\beta,\varepsilon)$ and vertices  
$u_1,\dots,u_h\in Y$ such that
\[
\left|\bigcup_{i=1}^{h} (N(u_i)\cap X)\right|>(1-\varepsilon)|X|.
\]
For each $i$, define
\[
X_i^{r}=\{x\in X: x u_i\in E(R)\},
\text{ and }
X_i^{b}=\{x\in X: x u_i\in E(B)\}.
\]
Then
\[
\left|\bigcup_{i=1}^{h}(X_i^{r}\cup X_i^{b})\right|
> (1-\varepsilon)|X|.
\]

Since $G$ contains no monochromatic triangle with two vertices in $X$ and one in $Y$,  
Observation~\ref{K3:free} implies that for all $i\neq j$,
\[
E(X_i^{r}) \cap E(R) = E(X_i^{b}) \cap E(B) = \varnothing, \text{ and }
 |X_i^{r}\cap X_j^{b}|<\alpha n.
\]
Moreover, also by Observation~\ref{K3:free},
\[
|N_R(y)\cap X_i^{b}|<\alpha n
\quad and \quad
|N_B(y)\cap X_i^{r}|<\alpha n
\]
for any $i \in [h]$ \text{and} $y\in Y$. Summing over $i$ gives
\[
|N_R(y)\cap X^{b}|<h\alpha n
\quad \text{and} \quad
|N_B(y)\cap X^{r}|<h\alpha n.
\]
Therefore, summing over $y \in Y$ deduces
\[
|E(X^{b},Y)\cap E(R)|<h\alpha n^{2}
\quad \text{and} \quad
|E(X^{r},Y)\cap E(B)|<h\alpha n^{2}.
\]
This completes the proof.
\end{proof}

Apply Claim~\ref{claim:two-part} separately to these two ordered pairs $(P,Q)$ and $(Q,P)$,
we obtain that there exist integer $h = h(\beta,\varepsilon)$ and
subsets
\[
X_1^{r}, X_2^{r}, \ldots, X_h^{r};\,
X_1^{b}, X_2^{b}, \ldots, X_h^{b} \subseteq X
\]
for $X \in \{P,Q\}$ such that all three items in Claim~\ref{claim:two-part} hold. 
In particular, (3) of Claim~\ref{claim:two-part} states that
\begin{align}\label{redbluebd:ieq}
|E(P^{r}, Q)\cap E(R)|,|E(P^{b}, Q)\cap E(B)|,|E(P, Q^{r})\cap E(R)|,|E(P, Q^{b})\cap E(B)| \le h\alpha n^2.
\end{align}
Note that since $h=h(\beta,\varepsilon)$ and recall the choice of parameters at the beginning of this proof, we have
\begin{align}\label{para:ieq}
    \alpha \ll \beta,\varepsilon,\frac mn, \frac1h.
\end{align}
Then
\begin{align}
    |E_G(P^r,Q^b)| &= |E_G(P^r,Q^b)\cap E(R)|+|E_G(P^r,Q^b)\cap E(B)|\nonumber\\
    &< |E_G(P^r,Q)\cap E(R)|+|E_G(P,Q^b)\cap E(B)|\nonumber\\
    &\le 2h\alpha n^2\nonumber\\
    &< \frac12\beta\varepsilon^2(\frac mn)^2n^2\nonumber\\
    &< \varepsilon(\beta-2\sqrt{\varepsilon})(\sqrt\varepsilon m)^2\nonumber\\
    &\le \varepsilon(\beta-2\sqrt{\varepsilon})|P||Q|\label{alpha:ieq}.
\end{align}

Recall that the pair $(P,Q)$ is $\sqrt{\varepsilon}$-regular, it is clear that either
$|P^r|\le \sqrt{\varepsilon}|P|$ or $|Q^b|\le\sqrt{\varepsilon}|Q|$; otherwise, we can deduce that the density $d_G(P^r,Q^b)\ge d_G(P,Q)-\sqrt{\varepsilon} > \beta-2\sqrt\varepsilon$.
Thus we have
\begin{align*}
    |E_G(P^r,Q^b)|=d_G(P^r,Q^b)\cdot|P^r||Q^b| > \varepsilon(\beta-2\sqrt\varepsilon) |P||Q|,
\end{align*}
which contradicts with (\ref{alpha:ieq}).
Similarly, we have $|E_G(P^b,Q^r)| \le 2h\alpha n^2$, and then either
$|P^b|\le \sqrt{\varepsilon}|P|$ or $|Q^r|\le\sqrt{\varepsilon}|Q|$.

According to (1) of Claim~\ref{claim:two-part}, we know that
$$
|P^r|+|P^b|\ge|P^r \cup P^b| > (1-\sqrt{\varepsilon})|P|
$$
and
$$
|Q^r|+|Q^b|\ge|Q^r \cup Q^b| \ge (1-\sqrt{\varepsilon})|Q|.
$$
This concludes that either
$|P^r| \le \sqrt{\varepsilon}|P|$ and $|Q^r| \le \sqrt{\varepsilon}|Q|$, or
$|P^b| \le \sqrt{\varepsilon}|P|$ and $|Q^b| \le \sqrt{\varepsilon}|Q|$.
Without loss of generality, we may assume that the first case holds. 
Since $|P^r\cup P^b|>(1-\sqrt{\varepsilon})|P|$ and $|Q^r\cup Q^b|>(1-\sqrt{\varepsilon})|Q|$, we have 
\begin{align}\label{AbBb:eq}
|P^b|\ge (1-2\sqrt{\varepsilon})|P| \quad \text{and} \quad |Q^b|\ge (1-2\sqrt{\varepsilon})|Q|.
\end{align}
Note that both $(P,S)$ and $(Q,S)$ are $\sqrt{\varepsilon}$-regular pairs with density at least $\beta-\varepsilon$. We apply Fact~\ref{deg:fact} to the pair $(P,S)$ with $Y=P^b$, this yields
\begin{align*}
    \lvert\{x \in S : \lvert N_G(x)\cap P^b\rvert \le (\beta-2\sqrt\varepsilon)|P^b|\}\rvert \le \sqrt\varepsilon |S|.
\end{align*}
Applying Fact~\ref{deg:fact} again to the pair $(Q,S)$ with $Y=Q^b$, we similarly obtain
\begin{align*}
    \lvert\{x \in S : \lvert N_G(x)\cap Q^b\rvert \le (\beta-2\sqrt\varepsilon)|Q^b|\}\rvert \le \sqrt\varepsilon |S|.
\end{align*}
Consequently, there exist at least $(1-2\sqrt{\varepsilon})|S|$ vertices $v \in S$ satisfying
\begin{align*}
 |N_G(v)\cap P^b|> (\beta-2\sqrt{\varepsilon})|P^b| \quad \text{and} \quad |N_G(v)\cap Q^b|> (\beta-2\sqrt{\varepsilon})|Q^b|.
\end{align*} 
Fix such a vertex $v \in S$.
We have the following claim.
\begin{claim}\label{nb(v)}
    $|N_B(v)\cap P^b|\ge \frac12\beta|P^b|$ and $|N_B(v)\cap Q^b|\ge \frac12\beta|Q^b|$.
\end{claim}
\begin{proof}
    We only prove the first inequality, and the second can be proved analogously. Assume that $|N_B(v)\cap P^b|< \frac12\beta|P^b|$ for contradiction, then we have
    $|N_R(v)\cap P^b|= |N_G(v)\cap P^b|-|N_B(v)\cap P^b|> (\beta-2\sqrt{\varepsilon})|P^b|-\frac12\beta|P^b| >\frac13\beta |P^b|$. Therefore, 
    $$
    \sum_{i=1}^h|N_R(v)\cap P_i^b| 
    \ge |N_R(v)\cap P^b|
    > \frac13\beta |P^b|.
    $$
    By the Pigeonhole Principle, there exists some $i\in [h]$ such that
    $$
    |N_R(v)\cap P_i^b|>\frac{\beta}{3h}|P^b| > \frac{\beta\sqrt\varepsilon}{3h}(1-2\sqrt{\varepsilon}) \frac mn n>\alpha n> \alpha(G),
    $$
    where the second inequality follows from (\ref{para:ieq}).
    But this leads to a contradiction as: Note that by Observation~\ref{K3:free},
    $$E_G(N_R(v)\cap P)\cap E(R)=\emptyset,$$
   and by (2) of Claim~\ref{claim:two-part},
    $$E_G(P_i^b)\cap E(B)=\emptyset.$$ 
Consequently, $E_G(N_R(v)\cap P_i^b)=\emptyset$. This contradicts with $|N_R(v)\cap P_i^b|>\alpha(G)$. 
Therefore, we have $|N_B(v)\cap P^b|= |N_G(v)\cap P^b|-|N_R(v)\cap P^b| \ge (\beta-2\sqrt{\varepsilon})|P^b|-\frac13\beta|P^b| \ge \frac12\beta|P^b|$.
\end{proof}

By Claim~\ref{nb(v)}, we have $|N_B(v)\cap P^b|\ge \frac12\beta|P^b| >\sqrt{\varepsilon}|P|$ and $|N_B(v)\cap Q^b|\ge \frac12\beta|Q^b|>\sqrt{\varepsilon}|Q|$. According to the regularity of $(P,Q)$, we have
$$
|E_G(N_B(v)\cap P^b,N_B(v)\cap Q^b)| \ge \frac{\beta^2}{4}(\beta-2\sqrt{\varepsilon})|P^b||Q^b|.
$$
On the other hand, 
$$
|E_G(N_B(v)\cap P^b,N_B(v)\cap Q^b)\cap E(R)| \le |E_G(P^b,Q)\cap E(R)|\le h\alpha n^2,
$$
where the last inequality follows from (3) of Claim~\ref{claim:two-part}.
Consequently, we have
\begin{align*}
    |E_G(N_B(v)\cap P^b,N_B(v)\cap Q^b)\cap E(B)|&\ge \frac{\beta^2}{4}(\beta-2\sqrt{\varepsilon})|P^b||Q^b|- h\alpha n^2\\
    &\ge \frac{\beta^2}{4}(\beta-2\sqrt{\varepsilon})(1-2\sqrt{\varepsilon})^2|P||Q|- h\alpha n^2\\
    &\ge \frac{\beta^2}{4}(\beta-2\sqrt{\varepsilon})(1-2\sqrt{\varepsilon})^2\varepsilon m^2- h\alpha n^2\\
    &>0,
\end{align*}
where the second inequality holds by \eqref{AbBb:eq}, and the last inequality holds by (\ref{para:ieq}).
Therefore, we can choose a blue edge
$$
uw\in E_G(N_B(v)\cap P^b,N_B(v)\cap Q^b)\cap E(B).
$$
Thus we obtain a blue triangle $uvw$ with $u,v, w$ located in the distinct sets $P, Q$ and $S$, respectively. But this is precisely configuration (5), which  was listed as impossible, yielding a contradiction.
\end{proof}

Now we are ready to give the proof of Lemma~\ref{lem:k3-tiling-in-F}.
\mainlem*

\begin{proof}
First, we consider the induced subgraph $H:=G[V_1\cup V_2\cup V_3]$. According to Lemma~\ref{lem:three-part}, we can take a maximum monochromatic $K_3$-tiling $\Gamma$ of $(H,\phi|_H)$ such that each $K_3$-copy intersects with $V_1$ in at most one vertex, while intersects with both $V_2$ and $V_3$ in at most two vertices.
Let $V_i'=V_i\setminus V(\Gamma)$ for $i=1,2,3$.
We claim that 
\begin{align}\label{size:reg:eq}
  \text{there exists some $i \in [3]$ such that }  |V_i'|<\sqrt{\varepsilon}m.
\end{align}
Otherwise, we apply Fact~\ref{slicing:fact} with $A=V_i$, $B=V_j$, $d=\beta$ and $\alpha=\sqrt{\varepsilon} > \varepsilon$ for any $1 \le i <j \le 3$. Note that $\varepsilon'=\max\{\varepsilon/\alpha,2\varepsilon\}=\sqrt{\varepsilon}$ since $\varepsilon \ll 1$. Hence, the pair $(V_i',V_j')$ is $\sqrt{\varepsilon}$-regular and its density satisfies $d(V_i',V_j') > \beta-\varepsilon$.

We can continue to apply Lemma~\ref{lem:three-part} to the graph $H':=H\setminus V(\Gamma)$ with partition $V(H')=V_1'\cup V_2' \cup V_3'$.
Hence there exists some monochromatic $K_3$ in $(H',\phi_{H'})$ who intersects with $V_1$ in at most one vertex, while intersects with both $V_2$ and $V_3$ in at most two vertices. This deduces that we can find a larger tiling by add this $K_3$ into $\Gamma$, contradicts with the maximum of $\Gamma$. 

If $|V_1\setminus V(\Gamma)| <\sqrt{\varepsilon}m$, then 
$|V_1\cap V(\Gamma)|\ge (1-\sqrt{\varepsilon})m,$ 
since each triangle of $\Gamma$ intersects $V_1$ in at most one vertex. As we desired.

Now we assume that $|V_1\setminus V(\Gamma)| \ge \sqrt{\varepsilon}m$. By \eqref{size:reg:eq} and the symmetry of $V_2$ and $V_3$, We can further assume that $|V_2\setminus V(\Gamma)| <\sqrt{\varepsilon}m$. Since each triangle of $\Gamma$ intersects $V_2$ in at most two vertices, 
$$
|\Gamma|\ge \frac12|V_2\cap V(\Gamma)|\ge \frac{1-\sqrt{\varepsilon}}2m.
$$ 
Let $V_1':=V_1\setminus V(\Gamma)$. Then $|V_1'|\ge  \sqrt{\varepsilon}m$. Consider the induced subgraph $Q:=G[V_1'\cup V_4\cup V_5]$. Take a maximum monochromatic $K_3$-tiling $\Gamma'$ in $(Q,\phi|_Q)$ such that each $K_3$-copy intersects with $V_1'$ in at most one vertex, while intersects with both $V_4$ and $V_5$ in at most two vertices.
With similar discussion as in the statement (\ref{size:reg:eq}),
at least one of
$$
|V_1'\setminus \Gamma'|,|V_4\setminus \Gamma'|,|V_5\setminus \Gamma'|
$$
have size less than $\sqrt{\varepsilon}m$. Since $\Gamma$ and $\Gamma'$ are disjoint, we can union them together to form a larger tiling $\Gamma\cup \Gamma'$.

If $|V_1'\setminus \Gamma'| < \sqrt{\varepsilon}m$, then, as 
$\Gamma\cup \Gamma'$ intersects with $V_1$ in at most one vertex, we have
$$
|\Gamma\cup \Gamma'|\ge  |V_1|-|V_1'\setminus \Gamma'| \ge (1-\sqrt{\varepsilon})m,
$$
as desired.

Now we assume that $|V_1'\setminus \Gamma'| \ge \sqrt{\varepsilon}m$ and further assume that $|V_4\setminus \Gamma'| < \sqrt{\varepsilon}m$ by the symmetry of $V_4$ and $V_5$.
Since each triangle in $\Gamma'$ intersects $V_4$ in at most two vertices, we have 
\begin{align*}
|\Gamma'|\ge \frac12|V_4\cap V(\Gamma)|\ge \frac{1-\sqrt{\varepsilon}}2m.
\end{align*}
Consequently, we have $|\Gamma\cup \Gamma'|\ge |\Gamma|+|\Gamma'|\ge (1-\sqrt{\varepsilon})m$, which completes the proof.

\end{proof}

\section{Concluding remarks}


In this paper, we present the Ramsey--Tur\'{a}n type theorem for monochromatic triangle-tilings.
It is natural to ask how far the present methods extend beyond triangles. For larger cliques, neighborhood intersections from different color classes interact in more delicate ways, and new ideas may be required even to force a single monochromatic $K_\ell$ under the Ramsey--Tur\'an condition. This leads to the following problem.

\begin{prob}
For each fixed $\ell \ge 4$, determine the optimal minimum-degree condition that guarantees the existence of a monochromatic $K_\ell$ or a monochromatic $K_\ell$-tiling of linear size in every sufficiently large $2$-edge-coloured graph $G$ with $\alpha(G)=o(n)$.
\end{prob}

\end{document}